\documentclass[11pt]{article}

\usepackage{amsfonts,amsmath,latexsym,color,epsfig}
\newtheorem{theorem}{Theorem}
\newtheorem{lemma}{Lemma}

\newenvironment{proof}
      {\medskip\noindent{\bf Proof:}\hspace{1mm}}
      {\hfill$\Box$\medskip}

\input{epsf}

\makeatletter
\def\Ddots{\mathinner{\mkern1mu\raise\p@
\vbox{\kern7\p@\hbox{.}}\mkern2mu
\raise4\p@\hbox{.}\mkern2mu\raise7\p@\hbox{.}\mkern1mu}}
\makeatother

\title{Large almost monochromatic subsets in hypergraphs}
\author{David Conlon\thanks{St John's College, Cambridge, United Kingdom.
E-mail: {\tt D.Conlon@dpmms.cam.ac.uk}. Research supported by a
research fellowship at St John's College.} \and Jacob
Fox\thanks{Department of Mathematics, Princeton, Princeton, NJ.
Email: {\tt jacobfox@math.princeton.edu}. Research supported by an
NSF Graduate Research Fellowship and a Princeton Centennial
Fellowship.} \and Benny Sudakov\thanks{Department of Mathematics,
UCLA,  Los Angeles, CA 90095. Email: {\tt bsudakov@math.ucla.edu}. Research supported in part by NSF CAREER award DMS-0812005 and by
USA-Israeli BSF grant.}}
\date{}
\begin{document}
\maketitle

\begin{abstract}
We show that for all $\ell$ and $\epsilon>0$ there is a constant $c=c(\ell,\epsilon)>0$ such that
every $\ell$-coloring of the triples of an $N$-element set contains a subset $S$
of size $c\sqrt{\log N}$ such that at least
$1-\epsilon$ fraction of the triples of $S$ have the same color. This result is tight up to the constant $c$ and answers an open question of Erd\H{o}s and Hajnal
from 1989 on discrepancy in hypergraphs. For $\ell \geq 4$ colors, it is
known that there is an $\ell$-coloring of the triples of an $N$-element set
whose largest monochromatic subset has cardinality only
$\Theta(\log \log N)$. Thus, our result demonstrates that the maximum almost monochromatic subset that an $\ell$-coloring
of the triples must contain is much larger than the corresponding monochromatic subset. This is in striking contrast with graphs, where these two
quantities have the same order of magnitude. To prove our result, we obtain a new upper bound on the $\ell$-color Ramsey
numbers of complete multipartite $3$-uniform hypergraphs, which answers another open question of Erd\H{o}s and Hajnal.
\end{abstract}

\section{Introduction}

The {\it Ramsey number} $r(n)$ is the smallest integer $N$ such that
every $2$-coloring of the edges of the complete graph on $N$ vertices contains a monochromatic clique of size $n$.
Ramsey's theorem states that $r(n)$ exists for all $n$. Determining or estimating Ramsey numbers is one of the
central problems in combinatorics, see the book Ramsey theory
\cite{GRS90} for details. A classical result of Erd\H{o}s and
Szekeres~\cite{ES35}, which is a quantitative version of Ramsey's
theorem, implies that $r(n) \leq 2^{2n}$ for every positive
integer $n$. Erd\H{o}s~\cite{E47} showed using probabilistic
arguments that $r(n) > 2^{n/2}$ for $n
> 2$. Over the last sixty years, there have been several
improvements on these bounds (see, e.g., \cite{C08}). However,
despite efforts by various researchers, the constant factors in
the above exponents remain the same.

Although already for graph Ramsey numbers there are significant
gaps between lower and upper bounds, our knowledge of hypergraph
Ramsey numbers is even weaker. The Ramsey number $r_k(n)$ is the
minimum $N$ such that every $2$-coloring of the $k$-tuples of an $N$-element set
contains a monochromatic set of size $n$, where a set is called monochromatic if all its $k$-tuples have the same color.
Erd\H{o}s, Hajnal, and Rado \cite{EHR65} showed that there are positive
constants $c$ and $c'$ such that $$2^{cn^2}<r_3(n)<2^{2^{c'n}}.$$
They also conjectured that $r_3(n)>2^{2^{cn}}$ for some constant
$c>0$ and Erd\H{o}s (see, e.g. \cite{CG98}) offered a \$500 reward for a proof. Similarly,
for $k \geq 4$, there is a difference of one exponential between
known upper and lower bounds for $r_k(n)$, i.e., $$t_{k-1}(cn^2)
\leq r_k(n) \leq t_k(c'n),$$ where the tower function $t_k(x)$ is
defined by $t_1(x)=x$ and $t_{i+1}(x)=2^{t_i(x)}$.

The study of $3$-uniform hypergraphs is particularly important for
our understanding of hypergraph Ramsey numbers. This is because of
an ingenious construction called the stepping-up lemma due to
Erd\H{o}s and Hajnal (see, e.g., Chapter 4.7 in \cite{GRS90}).
For $k \geq 3$, their method allows one to construct lower bound colorings for
uniformity $k+1$ from colorings for uniformity $k$, effectively
gaining an extra exponential each time it is applied. Therefore, proving that $r_3(n)$ has doubly exponential growth
will allow one to close the gap between the upper and lower bounds
for $r_k(n)$ for all uniformities $k$.

Despite the fact that Erd\H{o}s \cite{E90} (see also the book
\cite{CG98}) believed $r_3(n)$ is closer to $2^{2^{cn}}$, together
with Hajnal \cite{EH89}, he discovered the following interesting
fact about hypergraphs which maybe indicates the opposite. They
proved that there are $c,\epsilon>0$ such that every $2$-coloring
of the triples of an $N$-element set contains a subset $S$ of size $s> c(\log
N)^{1/2}$ such that at least $(1/2+\epsilon){s \choose 3}$ triples of $S$ have the same color. That is, this subset deviates from
having density $1/2$ in each color by at least some fixed positive
constant. Erd\H{o}s \cite{E94} further remarks that he would begin to doubt
that $r_3(n)$ is double-exponential in $n$ if one can prove that
any $2$-coloring of the triples of an $N$-set contains some set of
size $s=c(\epsilon)(\log N)^{\delta}$ for which at least $(1-\epsilon){s
\choose 3}$ triples have the same color, where $\delta>0$ is an absolute constant.
Erdos and Hajnal proposed \cite{EH89} that such a statement may even be true with $\delta =
1/2$. Our first result shows that this is indeed the case.

\begin{theorem}\label{maintheorem}
For each $\epsilon > 0$ and $\ell$, there is $c=c(\ell,\epsilon)>0$ such that
every $\ell$-coloring of the triples of an $N$-element set contains a
subset $S$ of size $s=c\sqrt{\log N}$ such that at least $(1-\epsilon){s \choose 3}$ triples of $S$ have the same color.
\end{theorem}

It is easy to see that this theorem is tight up to the constant factor $c$.
Indeed, consider an $\ell$-coloring of the triples of an $N$-element set in which every
triple gets one of $\ell$ colors uniformly at random. Using a standard tail estimate for the binomial distribution,
one can show that in this coloring, with high probability, every subset of size $\gg \sqrt{\log N}$ has a $1/\ell+o(1)$ fraction of its triples in each color.

Our result also shows a significant difference between the discrepancy problem in graphs and that in hypergraphs.
The {\it $\ell$-color} Ramsey number $r_k(n;\ell)$ is the minimum $N$ such that
every $\ell$-coloring of the $k$-tuples of an $N$-element
set contains a monochromatic set of size $n$. Erd\H{o}s and Hajnal
(see, e.g., \cite{GRS90}) constructed a $4$-coloring of the
triples of a set of size $2^{2^{cn}}$ which does not contain a
monochromatic subset of size $n$. In \cite{CFS08} we constructed a $3$-coloring of the
triples of a set of size $2^{n^{c\log n}}$ which does not contain a
monochromatic subset of size $n$. Thus, Theorem \ref{maintheorem} demonstrates (at least for $\ell \geq 3$)
that the maximum almost monochromatic subset that an $\ell$-coloring
of the triples must contain is much larger than the corresponding monochromatic subset. This is in striking contrast with graphs, where these two quantities have the
same order of magnitude, as demonstrated by a random $\ell$-coloring of the edges of a complete graph.

Another open problem from the 1989 paper of Erd\H{o}s and Hajnal \cite{EH89} asks whether one can exhibit a fixed hypergraph of density larger than
$1/2+\epsilon$ on $c\sqrt{\log N}$ vertices that occurs monochromatically.
That is, can we find dense hypergraphs with small Ramsey numbers?
We show that this is indeed the case by obtaining a new upper bound on the $\ell$-color Ramsey number of a complete multipartite $3$-uniform
hypergraph. A {\it hypergraph} $H=(V,E)$ consists of a vertex set $V$ and an edge set $E$, which is a collection of subsets of $V$. A hypergraph is {\it $k$-uniform} if each edge has exactly $k$ vertices. For a $k$-uniform hypergraph $H$, the Ramsey number $r(H;\ell)$ is
the minimum $N$ such that every $\ell$-coloring of the $k$-tuples of
an $N$-element set contains a monochromatic copy of $H$. The {\it complete $d$-partite} $k$-uniform hypergraph $K_d^k(n)$
is the $k$-uniform hypergraph whose vertex set consists of $d$
parts of size $n$ and whose edges are all $k$-tuples that have
their vertices in some $k$ different parts. The number of vertices of $K_d^3(n)$ is
$dn$ and the number of edges in $K_{d}^3(n)$ is ${d \choose 3}n^3
> (1-\frac{3}{d}){dn \choose 3}$, i.e., it has edge density more
than $1-\frac{3}{d}$. In particular, as $d$ increases, the edge
density of $K_d^3(n)$ tends to $1$. Therefore, Theorem
\ref{maintheorem} is an immediate corollary of the following
theorem.

\begin{theorem}\label{mainRamsey}
The $\ell$-color Ramsey number of the complete $d$-partite hypergraph $K_d^3(n)$ satisfies
$$r(K_{d}^3(n);\ell)  \leq 2^{\ell^{2r}n^2},$$ where $r=r_2(d-1;\ell)$ is the $\ell$-color Ramsey number of the complete graph on $d-1$ vertices.
\end{theorem}

In the next section, we discuss the technique of K\"ov\'ari, S\'os,
and Tur\'an on the classical problem of Zarankiewicz  and its
consequences, which will be helpful in proving our main theorems.
In Section \ref{proofsection}, we present the proof of Theorem
\ref{mainRamsey}. Throughout the paper, we systematically omit
floor and ceiling signs whenever they are not crucial for the sake
of clarity of presentation. We also do not make any serious
attempt to optimize absolute constants in our statements and
proofs. All logarithms in this paper are in base $2$.

\section{Complete bipartite graphs in dense graphs}

The problem of Zarankiewicz \cite{Z51} asks for the maximum number
$z(m,n;s,t)$ of edges in a bipartite graph $G$ which has $m$ vertices in
its first class, $n$ vertices in its second class, and does not contain a
complete bipartite graph $K_{s,t}$ with $s$ vertices in the first class
and $t$ in the other. In their celebrated paper, K\"ov\'ari, S\'os, and Tur\'an \cite{KST54}
used double counting together with the
pigeonhole principle to give a general upper bound on $z(m,n;s,t)$. Using this technique, we obtain the following
two simple lemmas which we need in the proof of our main result. The degree $d(v)$ of a vertex $v$ is the number of vertices
adjacent to $v$.

\begin{lemma}\label{mainlemma}

Let $G$ be a bipartite graph with parts $A$ and $B$ and with at least
$|A||B|/\ell$ edges. Then $G$ contains a complete bipartite subgraph with
one part having $a=|A|/\ell$ vertices from $A$ and the other part having
$b=2^{-|A|}|B|$ vertices from $B$.

\end{lemma}
\begin{proof}
Using the convexity of the function $f(x)={x \choose a}$ together with the fact that the average degree of a vertex in $B$ is at least
$a$, we conclude that the number of pairs $(A',v)$, where $A'$ is a subset of $A$ of size $a$ and $v
\in B$ is adjacent to every vertex in $A'$, is at least
$$\sum_{v \in B}{d(v) \choose a} \geq |B| {{\frac{1}{|B|}\sum_{v \in B} d(v)} \choose a}=|B|.$$
Since the set $A$ has at most $2^{|A|}$ subsets, by the pigeonhole principle, there are at least
$b=2^{-|A|}|B|$ vertices in $B$ adjacent to the same subset $A' \subset A$ of size $a$.
Together they form a complete bipartite graph with one part having $a$ vertices in $A$
and the other part having $b$ vertices in $B$, completing
the proof.
\end{proof}

\begin{lemma}\label{mainlemma1}
If a graph $G$ of order $n$ has $\epsilon n^2$ edges and $t
<\epsilon n$, then it contains $K_{s,t}$ with $s= \epsilon^t n$.
\end{lemma}
\begin{proof}
Note that the number of pairs $(U,v)$ with $U$ being a subset of $G$ of size $t$ and $v$ being a vertex of $G$ adjacent to every vertex in $U$ is
at least
$$\sum_{v \in V(G)}{d(v) \choose t} \geq n{2\epsilon n \choose t} \geq n(\epsilon n)^t/t!,$$
where we use convexity of $f(x)={x \choose t}$ and the fact that the average degree of a vertex in $G$ is $2\epsilon n >t+\epsilon n$.
If $G$ does not contain $K_{s,t}$ as a subgraph, then for every subset $U$ of $G$ of size $t$, there are at most $s-1$ vertices of $G$ adjacent to
$U$. Hence,
$$n(\epsilon n)^t/t! \leq (s-1){n \choose t}< sn^t/t! = n(\epsilon n)^t/t!,$$
a contradiction, which completes the proof.
\end{proof}

Both of these lemmas can be further improved
by using better upper bound estimates for binomial coefficients.
However, the above clean estimates are sufficient for our purposes.

\section{Proof of the main result} \label{proofsection}

First we briefly discuss the classical approach of Erd\H{o}s and
Rado which gives an upper bound on the Ramsey number of the complete $3$-uniform hypergraph on $d$ vertices.
Suppose the triples of a sufficiently large set $V$ of vertices
are $\ell$-colored. Let $r$ be the $\ell$-color Ramsey number of a complete graph of order $d-1$.
Erd\H{o}s and Rado greedily construct a set of
vertices $\{v_1,\ldots,v_{r+1}\}$ such that for any given
pair $1 \leq i <j \leq r$, all triples $\{v_i,v_j,v_k\}$
with $k>j$ are of the same color, which we denote by
$\chi(v_i,v_j)$. By definition of the Ramsey number $r=r_2(d-1;\ell)$, there is a
monochromatic clique of size $d-1$ in coloring $\chi$, and this
clique together with $v_{r+1}$ forms a monochromatic set of
size $d$ in the original coloring. Note that in their approach,
after having picked $\{v_1,\ldots,v_i\}$, we have a subset $S_i$
such that for any pair $a,b$ with $1 \leq a <b \leq i$, all
triples $\{v_a,v_b,w\}$ with $w \in S_i$ are the same color. The
subset $S_i$ consists of those vertices from which we can draw
future vertices.

Instead of picking vertices one by one as in the Erd\H{o}s-Rado
technique, we instead pick subsets one by one. After step $i$, we
will have disjoint subsets $V_{1,i},\ldots,V_{i,i}$ and another
subset $S_i$ such that for any pair $a,b$ with $1 \leq a
<b \leq i$, all triples in $V_{a,i} \times V_{b,i} \times V_{c,i}$
with $b<c \leq i$ and all triples in $V_{a,i} \times V_{b,i}
\times S_i$ are the same color. Similar to before, the subset $S_i$ consists of those vertices from which we can draw
future vertex subsets.

We now present the proof of Theorem \ref{mainRamsey}, which states that the $\ell$-color Ramsey number of the complete
$d$-partite hypergraph $K_d^3(n)$ satisfies $$r(K_{d}^3(n);\ell)  \leq 2^{\ell^{2r}n^2},$$ where $r=r_2(d-1;\ell)$. The
set of edges of a graph $G$ is denoted by $E(G)$ and $e(G)=|E(G)|$.

\noindent {\bf Proof of Theorem \ref{mainRamsey}:} Let
$r=r_2(d-1;\ell)$ and let $V$ be a set of
$N=2^{\ell^{2r}n^2}$ vertices whose triples are
$\ell$-colored. We will construct disjoint subsets
$V_1,V_2,\ldots,V_{r+1}$ of size $n$ such that for each pair $1
\leq i < j \leq r$ there is a color $\chi(i,j)$ for which all
triples in $V_i \times V_j \times V_k$ with $j < k \leq r+1$ have
color $\chi(i,j)$. Note that $\chi$ is an $\ell$-coloring of the
edges of the complete graph with vertex set $\{1,\ldots,r\}$ and therefore, from
the definition of the Ramsey number $r=r_2(d-1;\ell)$, it follows
that there are $d-1$ numbers $1 \leq i_1 < i_2 < \ldots < i_{d-1}
\leq r$ that make a monochromatic clique in coloring $\chi$. Then, using the properties of $\chi$
it is easy to see that
the corresponding sets $V_{i_1},V_{i_2},\ldots,V_{i_{d-1}}$ together with $V_{r+1}$
make a monochromatic $K_d^{3}(n)$. So we are left with constructing subsets
$V_1,\ldots,V_{r+1}$ with the desired properties, which we will do in $r$ rounds.

In the first round we pick $V_{1,1} \subset V$ of size
$\ell^{-1}\sqrt{\log N}$ arbitrarily and let $S_1=V \setminus
V_{1,1}$, so $|S_1| = N-\ell^{-1}\sqrt{\log N} \geq N^{3/4}$.
Now suppose that we already constructed disjoint subsets
$V_{1,i},V_{2,i},\ldots,V_{i,i}$ each of size $\ell^{-i}\sqrt{\log
N}$ and a subset $S_i$ of vertices disjoint from
$V_{1,i},\ldots,V_{i,i}$ with $|S_i| \geq N^{1/4+2^{-i}}$ such
that if $1 \leq a < b \leq i$, then all triples in $V_{a,i} \times
V_{b,i} \times V_{c,i}$ with $b < c \leq i$ and all triples in
$V_{a,i} \times V_{b,i} \times S_i$ have color $\chi(a,b)$. Note
that this is satisfied for $i=1$.

We next show how we proceed through round
$i+1$. Consider all triples with one vertex in
$V_{1,i}$ and the other two vertices in $S_i$.
By the pigeonhole principle, at least a $1/\ell$ fraction of these triples
have the same color which we denote by $\chi(1,i+1)$. Let $H$ be an auxiliary
bipartite graph whose first part $A$ is $V_{1,i}$, second part $B$
consists of all the unordered pairs $(w, w')$ from $S_i$, and whose edges are those
pairs $(a,b), a \in A$ and $b=(w, w') \in B$ such that the triple $(a, w, w')$ has color
$\chi(1,i+1)$. By applying Lemma
\ref{mainlemma} to $H$ we find a subset $V_{1,i+1} \subset V_{1,i}$
with $|V_{1,i+1}|=|V_{1,i}|/\ell=\ell^{-(i+1)}\sqrt{\log N}$
and a graph $G_{1,i}$ on $S_i$ with at least
$$2^{-|V_{1,i}|}{|S_i| \choose 2} \geq
2^{-2-\ell^{-i}\sqrt{\log N}}|S_i|^2$$ edges such that all
triples consisting of a vertex from $V_{1,i+1}$ and an edge from
$G_{1,i}$ have color $\chi(1,i+1)$.

Continuing this process, suppose that
after $j$ steps, we have already picked for all $1 \leq h \leq j$
the sets $V_{h,i+1} \subset V_{h,i}$ with $|V_{h,i+1}| =
|V_{h,i}|/\ell=\ell^{-(i+1)}\sqrt{\log N}$ and graphs $G_{j,i} \subset G_{j-1,i}
\subset \ldots \subset G_{1,i}$ on $S_i$ which have the following properties. The number of
edges of $G_{j,i}$ is at least $$2^{-2-j\ell^{-i}\sqrt{\log
N}}|S_i|^2$$ and every triple consisting of a vertex in
$V_{h,i+1}$ together with the two vertices of an edge from
$G_{h,i}$ have color $\chi(h,i+1)$.

In the step $j+1$, consider all triples whose one vertex is in
$V_{j+1,i}$ and the other two vertices form an edge of $G_{j,i}$.
By the pigeonhole principle, at least a $1/\ell$ fraction of these triples
have the same color, which we denote by $\chi(j+1,i+1)$. Similar to before, let $H$ be  an auxiliary
bipartite graph whose first part $A$ is $V_{j+1,i}$, second part $B$
consists of all the edges of $G_{j,i}$, and the edges of $H$ are those
pairs $(a,b), a \in A$ and $b=(w, w') \in E(G_{j,i})$ such that the triple $(a, w, w')$ has color
$\chi(j+1,i+1)$. By applying Lemma
\ref{mainlemma} to $H$ we find a subset
$V_{j+1,i+1} \subset V_{j+1,i}$ with
$|V_{j+1,i+1}|=|V_{j+1,i}|/\ell=\ell^{-(i+1)}\sqrt{\log N}$
 and a subgraph $G_{j+1,i} \subset G_{j,i} \subset \ldots \subset
G_{1,i}$ on $S_i$ with at least $$2^{-|V_{j+1,i}|}e(G_{j,i})
\geq 2^{-|V_{j+1,i}|}2^{-2-j\ell^{-i}\sqrt{\log N}}|S_i|^2 =
2^{-2-(j+1)\ell^{-i}\sqrt{\log N}}|S_i|^2$$ edges such that all
triples whose one vertex is from $V_{j+1,i+1}$ and the other two
vertices form an edge from $G_{j+1,i}$ have color $\chi(j+1,i+1)$.

After $i$ such steps, we have $V_{1,i+1},\ldots,V_{i,i+1}$ each of
size $\ell^{-(i+1)}\sqrt{\log N}$ and a sequence of graphs
$G_{i,i} \subset G_{i-1,i} \subset \ldots \subset G_{1,i}$ on
$S_i$ such that the number of edges of $G_{i,i}$ is at least
$$2^{-2-i\ell^{-i}\sqrt{\log N}}|S_i|^2 \geq 2^{-\sqrt{\log N}}|S_i|^2,$$ and, for $1 \leq h \leq i$,
all triples consisting of a vertex in $V_{h,i+1}$ together with an
edge from $G_{h,i}$ are color $\chi(h,i+1)$.

Now we apply Lemma \ref{mainlemma1} (with $\epsilon=2^{-\sqrt{\log N}}$ and $t=\ell^{-(i+1)}\sqrt{\log N}$) to graph $G_{i,i}$ to find disjoint subsets
$V_{i+1,i+1}$ and $S_{i+1}$ of $S_i$ that form a complete
bipartite graph in $G_{i,i}$ and satisfy
$|V_{i+1,i+1}|=\ell^{-(i+1)}\sqrt{\log N}$, and
\begin{eqnarray*} |S_{i+1}| & \geq &  \left(2^{-\sqrt{\log
N}}\right)^{|V_{i+1,i+1}|}|S_i| = N^{-\ell^{-(i+1)}}|S_i| \geq
N^{-2^{-(i+1)}}N^{1/4+2^{-i}} \geq
N^{1/4+2^{-(i+1)}}.\end{eqnarray*}
By construction we have that for $j<i+1$ all triples in
$V_{j,i+1}\times V_{i+1,i+1}\times S_{i+1}$ have color $\chi(j,i+1)$.
On the other hand, since $V_{j,i+1} \subset V_{j,i}$ for $j<i+1$ and
$V_{i+1,i+1},  S_{i+1} \subset S_i$ we have by induction that if
$1 \leq a < b \leq i$, then all triples in $V_{a,i+1} \times
V_{b,i+1} \times V_{c,i+1}$ with $b < c \leq i+1$ and all triples in
$V_{a,i+1} \times V_{b,i+1} \times S_{i+1}$ have color $\chi(a,b)$.
This completes the description of round $i+1$ of our induction process.

After $r$
such iterations, we have
disjoint subsets $V_{1,r},\ldots,V_{r,r}$, each of size
$\ell^{-r}\sqrt{\log N}=n$ and $S_{r}$ of size at least $N^{1/4}
\geq n$ such that if $1 \leq a < b \leq r$, then all triples in
$V_{a,r} \times V_{b,r} \times V_{c,r}$ with $b < c \leq r$ and
all triples in $V_{a,r} \times V_{b,r} \times S_r$ have color
$\chi(a,b)$. Letting $V_{i}=V_{i,r}$ for $1 \leq i \leq r$ and
$V_{r+1}$ to be any subset of $S_r$ of size $n$, we obtain the sets
$V_1,\ldots,V_{r+1}$ with desired properties. This completes the proof. \hfill$\Box$\medskip

\section{Concluding remarks}
\begin{itemize}
\item
It would be very interesting to extend Theorems \ref{maintheorem} and \ref{mainRamsey} to uniformity $k \geq 4$. In
\cite{CFS08} we obtain some preliminary remarks in this direction. We show that for all $k,\ell$ and $\epsilon>0$ there is
$\delta=\delta(k,\ell,\epsilon)>0$ such that every $\ell$-coloring of the $k$-tuples of an $N$-element set contains a subset of size $s=(\log N)^{\delta}$ which
contains at
least $(1-\epsilon){s \choose k}$ $k$-tuples of the same color.  Unfortunately, notice that $\delta$ here depends on $\epsilon$. Just as we deduce Theorem \ref{maintheorem} from Theorem \ref{mainRamsey}, this result is obtained by proving an upper bound on the Ramsey numbers of complete multipartite $k$-uniform hypergraphs.

\item It would be nice to determine the best possible dependence of $c$ on $\epsilon$ in Theorem \ref{maintheorem}.
As we already mention in the introduction, this theorem follows from our bound on the $\ell$-color Ramsey number of the complete
$d$-partite $3$-uniform hypergraph with $d=\Theta(\epsilon^{-1})$. From Theorem \ref{mainRamsey} we have that
$c \leq  \ell^{-r}$, where $r$ is the $\ell$-color Ramsey number of a complete graph of order $d-1$.
Therefore, using the simple upper bound $r=r_2(d-1;\ell) \leq \ell^{(d-1)\ell}$ we obtain that
$c(\ell,\epsilon)\leq 2^{-\ell^{\Theta(\ell/\epsilon)}}$. It seems likely that this double exponential dependence on $1/\epsilon$ is not correct.

\end{itemize}

\end{document}